\documentclass{article}
\usepackage{amsmath, amscd, amsthm}
\begin{document}
\title{Surgery on Foliations}
\author{Oliver Attie and Sylvain Cappell}
\maketitle
\begin{abstract}
  In this paper, we set up two surgery theories and two kinds of Whitehead
  torsion for foliations. First, we construct a bounded surgery theory and
  bounded Whitehead torsion for foliations, which correspond to the Connes'
  foliation algebra in the K-theory of operator algebras, in the sense that
  there is an analogy between surgery theory and index theory,  and
  a Novikov Conjecture for bounded surgery on foliations
  in analogy with the foliated Novikov conjecture of P.Baum and
  A.Connes in operator theory. This surgery theory classifies the leaves
  topologically. Secondly, we construct a bounded geometry
  surgery for foliations,
  which is a generalization of blocked surgery, and a bounded geometry
  Whitehead torsion. The classifications in this surgery theory
  include the specification of  the Riemannian metrics of
  the leaves up to quasi=isometry. We state Borel conjectures for
  foliations, which solves a
  problem posed by S.Weinberger \cite{Wein}, and verify these in some cases of
  geometrical interest. 
\end{abstract}
\section{Introduction}
In 1980 M.Pimsner and D.Voiculescu \cite{PV} proved that the K-theory of the 
$C^*$-algebra of the Kronecker foliation, the irrational 
rotation algebra, is isomorphic to the group $\textbf{Z}\oplus\textbf{Z}$ in both 
dimension 0 and 1. Pimsner and Voiculescu used their six-term exact sequence 
to prove this result. In this paper we develop a general surgery theory for foliations
and we prove the analogous result for this theory. We have two surgery
theories: one a bounded geometry foliated surgery theory, the
classifications of which include the specification of the Riemannian
metrics of the leaves up to quasi-isometry, and the other a boundedly controlled
foliated surgery theory where only the topology of the leaves is classified.
We prove first that $M \times \mathcal{F}$, where 
$\mathcal{F}$ is the Kronecker foliation and $M$ is a manifold of dimension
at least 5, has structure set $H_*(T^2;\mathcal{S}_*(M))$ where $T^2$ is 
the 2-torus, and $\mathcal{S}_*(M)$ the structure set of $M$, 
using codimension 1 splitting and the first author's results with S.Hurder \cite{AH}. We then 
introduce foliated Whitehead torsion, a foliated s-cobordism theorem
and surgery groups, analogous to Connes' transverse index theory \cite{Connes},
for foliations of arbitrary codimension, with the
constraint that the codimension is bounded below by the dimension of
the manifold so that the leaves are dimension $\ge 5$.
We prove a surgery exact sequence and state Novikov and Borel conjectures for 
foliations. 
P.Baum and A.Connes have stated a Novikov conjecture for foliations 
\cite{BC} using operator K-theory of the foliation algebra
which is analogous to our surgery theoretic version of the Novikov
Conjecture in the bounded case but is different
for the surgery theoretic case case where leaves are considered as
manifolds of bounded geometry. This surgery theory is a generalization
of the blocked surgery of \cite{Q}.
The two surgery theories are related by the fact that one gets the
bounded surgery theory by forgetting the Cheeger's finiteness condition
(that there are a finite number of types in a ball of fixed radius $r$,
for any $r$) along the leaves.
The theory of Baum and Connes using K-theory of operator algebras
corresponds to this bounded surgery theory.
The Novikov Conjecture of Baum and Connes involves the fundamental groupoid
of the foliation, whereas Novikov conjecture for boundedly controlled
foliated surgery involves the holonomy groupoid and the
locally finite $L$-homology of the leaves. These two conjectures
agree when the holonomy groupoid is the same as the 
fundamental groupoid of the leaves of the foliation.
We use bounded  surgery theory to prove this
conjecture in the case where $(M,\mathcal{F})$ is ultra-spherical.
S.Hurder \cite{Hurder} has proven the foliated Novikov conjecture for
a large class of foliations using his exotic index theory.

The bounded geometry surgery theory works is analogous to blocked
surgery, and we use this blocked surgery theory to study the problem of
when a manifold is a leaf of a foliation.
We note that, in the case of the Kronecker foliation,
although there are a continuously infinite number of quasi-isometry types 
of a leaf, recurrence of leaves makes the number of types of a foliation 
countable \cite{AH}. We see how this observation agrees with bounded
geometry surgery of the Kronecker foliation.

\section{Pimsner-Voiculescu for Structure Sets}
\newtheorem{dfn}{Definition}[section]
Throughout the paper we will work in the PL category. Our constructions can
be smoothed as well (See \cite{Att}). In addition we will assume that all
foliations have leaves of dimension $\ge 5$.
\begin{dfn} A foliated map is a map 
$$f:(M, \mathcal{F}) \to (N,\mathcal{G})$$
so that if $x$ and $y$ are on the same leaf, $f(x)$ and $f(y)$ are on the same 
leaf, and $f$ is a $bg$ map when restricted to each leaf.
A foliated homotopy is a foliated map which is a 
homotopy $f_t$ for which for each $t$, $f_0(x)$ is 
on the same leaf as $f_t(x)$. A foliated homotopy equivalence is one for
which $f \circ g$ is foliated homotopic to the identity and $g \circ f$ is 
foliated homotopic to the identity, where
$$g:(N,\mathcal{G}) \to (M,\mathcal{F})$$
is a foliated map.
\end{dfn}
\begin{dfn} A foliated simple homotopy equivalence
$$f:(M,\mathcal{F}) \to (N, \mathcal{G})$$
between two codimension one foliations is a foliated map which is a $bg$ 
simple homotopy equivalence \cite{Att} along the leaves and along the 
transverse foliation. A foliated simple homotopy equivalence
$$f:(M,\mathcal{F}) \to (N, \mathcal{G})$$
between two foliations of codimension $\ge 2$ is a foliated map which is a
$bg$ simple homotopy equivalence \cite{Att}
along the leaves and along the transversals of the foliation.
\end{dfn}
\begin{dfn} A foliated s-cobordism is a foliated cobordism 
$(W,\mathcal{H})$ between two foliations $(M,\mathcal{F})$ and 
$(N,\mathcal{G})$ so that $(W,\mathcal{H})$ is foliated simple homotopy 
equivalent to $(M,\mathcal{F})$ and to $(N,\mathcal{G})$. 
A foliated h-cobordism $(W,\mathcal{H})$ between two foliations
$(M,\mathcal{F})$ and $(N,\mathcal{G})$ so that $(W,\mathcal{H})$ is foliated
homotopy equivalent to either end.
\end{dfn}
\begin{dfn} Let $(M,\mathcal{F})$ be a foliation. 
The structure set of $(M,\mathcal{F})$ is the set of foliated simple 
homotopy equivalences $f:(M. \mathcal{F}) \to (N, \mathcal{G})$ with 
two foliated simple homotopy equivalences being the same if there is a 
foliated s-cobordism $W$ between $M$ and $N$.
\end{dfn}
\newtheorem{thm}{Theorem}[section]
\begin{thm}[Pimsner-Voiculescu for Structure Sets] The structure set of
the product of the Kronecker foliation $(T^2, \mathcal{F})$
 with a compact manifold $M^n$, $n \ge 5$ is 
$$\mathcal{S}^{foliated}_*(M \times \mathcal{F})=H_*(T^2;\mathcal{S}_*(M))$$
\end{thm}
\textit{Proof:} We recall the result of Attie-Hurder \cite{AH} that leaves of 
a codimension one foliation of  the form $M \times \textbf{R}$ are 
end-periodic.
We wish to show that leaves of a foliation of the
foliated simple homotopy type of the foliation $M \times \mathcal{F}$ are 
of the form $M \times \textbf{R}$. We observe that since the leaves are 
bg simple homotopy equivalent to $M \times \textbf{R}$ they are boundedly
controlled simple homotopy equivalent to $M \times \textbf{R}$. 
We can then apply boundedly controlled surgery and find, by applying the 
boundedly controlled Borel conjecture \cite{FP} for $\textbf{R}$ that the 
boundedly controlled structure set $\mathcal{S}^{bdd}(M \times \textbf{R})$ is 
isomorphic to $\mathcal{S}(M)$, which proves that the leaf is boundedly
controlled homeomorphic to $M^\prime \times \textbf{R}$ for some 
structure $M^\prime$ on $M$.  Hence the leaves are end-periodic. To study the
transverse direction we apply codimension 1 splitting to the 
transverse foliation, since we have a bg simple homotopy equivalence and 
the splitting obstruction vanishes and observe that by Attie-Hurder \cite{AH} 
we have end-periodicity in this direction as well. The result follows.
\section{Surgery on Foliations}
\begin{dfn}\cite{Haefliger} $T$ is a manifold of dimension equal to the codimension
  of the leaves, equipped with an immersion $j:T \to X$ transverse to the
  leaves and where the image meets each leaf at least once. The elements of
  $\Gamma_T$ are represents by the triples $(x,c,y)$ where $c$ is a
  path connected $j(y)$ to $j(x)$ in a leaf, two paths being equivalent
  if they determine the same holonomy. In other words $\Gamma_T$ is the
  subspace of $T \times G \times T$ formed from triples $(x,y,g)$ so
  that $j(x)=\beta(g)$ and $j(y)=\alpha(g)$.

  The foliation $\mathcal{F}$ can be defined by an open covering
  $\{U_i\}$ and submersions $f_i:U_i \to T_i$, where $T_i$ is a manifold
  of dimension equal to the codimension of the leaves, the $f_i$ being
  surjective with connected fibers. The $f_i$ should satisfy the following
  compatibility condition: for every pair $(i,j)$, there exists a
  homeomorphism $g_{ij}:f_j(U_i\cap U_j) \to f_i(U_i\cap U_j)$ so that
  $f_i=g_{ij} \cdot f_j$ on $U_i \cap U_j$. Let $T$ be the disjoint of
  $T_i$. The local homeomorphisms $g_{ij}$ generate a pseudogroup of
  $\mathcal{F}$ associated to the cocycle of definition $(U_i, f_i, g_{ij})$.
  We denote by $\Gamma_T$ the topological groupoid of germs of its
  elements.

  It is clear that we have a bijection between leaves of $\mathcal{F}$ and
  orbits of the groupoid $\Gamma_T$.
  To the transverse holonomy groupoid $\Gamma_T$ of a foliation $\mathcal{F}$
  on a manifold $X$ is associated, its classifying space $B\Gamma_T$ and
  a continuous map $X \to B\Gamma_T$.
\end{dfn}

\begin{thm}[s-cobordism theorem] A foliated s-cobordism $W$ where the 
leaves are of dimension $\ge$ 6 is foliated PL homeomorphic to a product.
\end{thm}
\textit{Proof:} Using the $bg$ s-cobordism theorem \cite{Att}, we see that the
leaves are products. To deal with the transversals we use the transversal
groupoid $\Gamma_T$ in \cite{Haefliger}, to give the transversal groupoid
control on the whole cobordism, via the product structure on the leaves.
We have that the foliation is defined by a map $W \to B\Gamma_T$ rel one end
and since the leaves are products, the transversals are products by the fact
that the map is unique up to homotopy. 
\par\noindent
\begin{dfn} The foliated Whitehead group of $(M, \mathcal{F})$ is 
the group of foliated h-cobordisms of $(M,\mathcal{F})$ and is denoted
$Wh_{\mathcal{F}}(M)$.
\end{dfn}
The following definition follows the definition of a foliated space
in \cite{CC}.
\begin{dfn} A foliated Poincare duality complex is a metric space
 $X$ which is a foliated space whose leaves are 
$bg$ Poincare duality complexes in the following sense. The foliation
$\mathcal{F}$ is an equivalence relation on $X$, equivalence 
classes being connected, embedded $bg$ Poincare duality spaces all of the same
dimension $k$. A foliated chart 
in $X$ is a pair $(U, \phi)$ where $U \subseteq X$ is open and 
$\phi:U \to T \times Z$ is a homeomorphism where $Z$ is an open neighborhood 
in a fixed metric space $Z_1$ and $T$ is an open ball in a given $bg$ 
Poincare duality space of finite radius. The set $P_y=\phi^{-1}(T \times y)$
where $y \in Z$ is called a plaque. If $P$ and $Q$ are plaques, then $P\cap Q$
is open in a given $bg$ Poincare duality space which is the same for $P$ and 
$Q$. The union of all plaques which are open in a given $bg$ Poincare duality 
space is that space itself.   
\end{dfn}
The following definitions are based on \cite{Wall} Chapters 9 and 10.
\par\noindent
\begin{dfn} Let $(M,\mathcal{F})$ be a foliation. The foliated
surgery group of $(M,\mathcal{F})$ is the group of foliated cobordism
classes of unrestricted objects where an unrestricted object is:
\begin{itemize}
\item A foliated Poincare pair $(Y,X)$ over $(M,\mathcal{F})$.
\item A foliated map $\phi:((W,\mathcal{G}), (\partial W, \mathcal{H})) \to
(Y,X)$ of pairs of degree 1, where $(W,\mathcal{G})$ is a foliation and 
$\phi\mid \partial W: \partial W \to X$ is a foliated simple homotopy 
equivalence.
\item A $bg$ stable framing $F$ of $\tau_W \oplus \phi^*(\tau)$, where 
$\tau$ is the Spivak normal fibration of $(Y,X)$.
\item A map $\omega : Y \to K$, where $K$ is a foliated complex so that
the pullback of the double cover of $K$ to $Y$ is orientation preserving.
\end{itemize}
We write $\theta \sim 0$ to denote that we can construct the triples:
\begin{itemize}
\item
a simple foliated Poincare triad $(Z;Y,Y_+)$ with $Y \cap Y_+=X$ and a bundle
$\mu$ over $Z$ extending $\nu$. 
\item a compact foliated manifold triad $(P;N,N_+)$ with $N \cap N_+=M$.
\item a map $\psi:(P;N, N_+) \to (Z;Y,Y_+)$ of degree 1 extending $\phi$
and inducing a foliated simple homotopy equivalence $N_+ \to Y_+$.
\end{itemize}
\par\noindent
The set of quadruples under the equivalence relation $\sim$ is the 
foliated surgery group.
We denote the foliated surgery group by $L_*^\mathcal{F}(M)$.
\end{dfn}
\begin{dfn} Let $X, \mathcal{F}$ be a foliated Poincare duality 
complex. The group of normal invariants of $X$, $NI_{\mathcal{F}}(X)$ is 
the bordism group of quadruples $(M,\phi,\nu,F)$. where $M$ is a foliation,
$\phi:M \to X$ of degree 1, $\nu$ a vector bundle over $X$, and trivialization
$F$ of $\tau_M \oplus \phi^*\nu$.
\end{dfn}
The formulation of the following conjectures is a problem in \cite{Wein}.
\par\noindent
\begin{dfn}[Borel Conjecture for Foliations] Let $(M, \mathcal{F})$ be a 
foliation. If the leaves of $\mathcal{F}$
are aspherical, then there is an isomorphism $\theta$
$$NI_{\mathcal{F}}(M) \to L_*^\mathcal{F}(M)$$
given by taking the surgery obstruction of a normal map in 
$NI_\mathcal{F}(M)$.
\end{dfn}
\begin{dfn}[Novikov Conjecture for Foliations] Let $(M,\mathcal{F})$ be a 
foliation. If the leaves of $\mathcal{F}$ are aspherical, then the above
map $\theta$ is injective. 
\end{dfn}
\begin{dfn}[Surgery Exact Sequence for foliations] The following sequence is 
exact:
$$... \to L_{*+1}^\mathcal{F}(M)\to \mathcal{S}_{\mathcal{F}}(M)\to 
NI_{\mathcal{F}} \to L_*^\mathcal{F}(M) \to ...$$
\end{dfn}
\section{Algebraic K-Theory}
\begin{dfn} Let $(M,\mathcal{F})$ be a manifold equipped with a foliation.
  The holonomy groupoid $H=Hol(M,\mathcal{F})$ is a smooth groupoid with
  $H_0=M$ as the space of objects. If $x,y \in M$ are two points on different
  leaves there are no arrows between $x$ and $y$ in $H$. If $x$ and $y$ lie
  on the same leaf $L$, an arrow $h:x \to y$ in $H$ is an equivalence class
  $h=[\alpha]$ of smooth paths $\alpha:[0,1] \to L$ with $\alpha(0)=x$ and
  $\alpha(1)=y$. To explain the equivalence relation, let $T_x$ and $T_y$
  be small $q$-disks through $x$ and $y$ transverse to the leaves of the
  foliation. If $x^\prime \in T_x$ is a point sufficiently close to $x$ on
  a leaf $L^\prime$, then $\alpha$ can be ``copied'' inside $L^\prime$ to
  give a path $\alpha^\prime$ near $\alpha$ with endpoint $y^\prime \in T_y$.
  In this way one obtains the germ of a diffeomorphism from $T_x$ to $T_y$
  sending $x$ to $y$ and $x^\prime$ to $y^\prime$. This germ is called the
  holonomy of $\alpha$ and denoted $hol(\alpha)$. Two paths $\alpha$ and
  $\beta$ from $x$ to $y$ in $L$ are equivalent, i.e. define the same
  arrow $x \to y$, if and only if $hol(\alpha)=hol(\beta)$. We obtain a well
  defined smooth groupoid $H=Hol(M,\mathcal{F})$. This groupoid is a foliation
  groupoid, and the (discrete) isotropy group $H_x$ at $x$ is called the
  holonomy group of the leaf through $x$.
  \end{dfn}
\begin{dfn} \cite{AM} Let $X_1$ and $X_2$ be spaces equipped with continuous maps
  $p_1, p_2$ to a metric space $Z$. Then a map $f:X_1 \to X_2$ is boundedly
  controlled if there exists an integer $m \ge 0$ so that for all $z \in Z$ of
  $r \ge 0$, $p_1^{-1}(B_r(z)) \subseteq f(p_2^{-1}(B_{r+m}(z)))$, where
  $B_r(z)$ denotes the metric ball in $Z$ of
  radius $r$ about $z$. Or, equivalently, there is a constant $m \ge 0$ so that
  $$dist_Z(p_2 \circ f(x), p_1(x))<m$$
  for all $x \in X_1$.
\end{dfn}
\begin{dfn} (See \cite{AM}). Let $(X_1,\mathcal{F}_1)$ and $(X_2,\mathcal{F}_2)$ be foliated
  spaces, and let $(Z,\mathcal{G})$ be a foliated space such that there are
  foliated maps $p_1:(X_1, \mathcal{F}_1)\to (Z,\mathcal{G})$ and
  $p_2:(X_2,\mathcal{F}_2)\to (Z,\mathcal{G})$. Then a foliated
  map $f:(X_1,\mathcal{F}_1)\to (X_2,\mathcal{F}_2)$ is boundedly controlled
  over the foliated space $(Z,\mathcal{G})$ such that for any leaf $L$ and
  for any transversal $T$ in $X_1$ there is a constant $m \ge 0$ so that
  $$dist_Z(p_2 \circ f(x), p_1(x))<m$$
  for $x \in L$ or $x \in T$.
  \end{dfn}
\begin{dfn} (See \cite{AM}). Let $(X,\mathcal{F})$ be a foliated space controlled over a
  metric space $Z$ by a map $p$. Denote by $\mathcal{P}$ the category of
  metric balls in $Z$ with morphisms given by inclusions. Define
  $\mathcal{P}Hol(X,\mathcal{F})$ to be the category whose objects are pairs
  $(x,K)$ where $K \in \mid \mathcal{P} \mid$ is an object of $\mathcal{P}$
  and if $x$ and $y$ are on the same leaf, a morphism $(x,K)\to (y,L)$ is a
  pair $(\omega, i)$ where $i \in \mathcal{P}(K,L)$ is a morphism in
  $\mathcal{P}$ from $K$ to $L$ and $\omega$ is a holonomy class of paths in
  $p^{-1}L$ from $y$ to $p^{-1}(i(x))$.
  \end{dfn}
\begin{dfn} (See \cite{AM}). Let $(S,\sigma)$ be a controlled basis. A controlled free
  $R\mathcal{P}Hol(X,\mathcal{F})$-module is defined so that given a
  controlled basis $(S,\sigma)$,
  \par\noindent i) For any metric ball $B$, $F(\sigma)(B)$ is the free
  $R$-module on the basis
  $$\{(\beta, s) \mid s\in S, \beta\mbox{ a path from }\sigma(s)\mbox{ to }b\}$$
  \par\noindent ii) The controlled holonomy group $\mathcal{P}Hol(M,\mathcal{F})$
  acts by composition of paths.
  \end{dfn}
The following bounded category corresponds to S.Hurder's exotic index theory
for foliations \cite{Hurder}.
\begin{dfn} The category of controlled free $\textbf{Z}\mathcal{P}Hol(X,
  \mathcal{F})$-modules is defined to be the category of controlled modules
  over the leaves of $\mathcal{F}$, with morphisms equal to the controlled
  morphisms. We denote this category by $\textbf{Z}\mathcal{P}Hol(X,
  \mathcal{F})^{bdd}$.
\end{dfn}
The $bg$ modules are new and would correspond to a foliated version of
Roe's uniform index theorem.
\begin{dfn} The category of controlled free $bg$
  $\textbf{Z}\mathcal{P}Hol(X,\mathcal{F})$-modules is defined to be the
  category of controlled modules so that in any ball of fixed radius of a leaf
  or transverse leaf the modules fall into a finite number of types. Morphisms
  are defined to be morphisms of the modules so that if the control space is
  partitioned into neighborhoods of a fixed radius, the restrictions fall into
  a finite number of equivalence classes. We will denote this category
  by $\textbf{Z}\mathcal{P}Hol(X,\mathcal{F})^{bg}$.
\end{dfn}
\begin{dfn} We define $K_1(\textbf{Z}\mathcal{P}Hol(X,\mathcal{F})^{bg})$ to
  be the abelian group generated by $[F,\alpha]$ where $F$ is a controlled free
  $bg$ module and $\alpha$ is an automorphism of $F$ so that
  \par\noindent (i) $[F,\alpha]=[F^\prime, \alpha^\prime]$ if there is an
  isomorphism $\phi:F\to F^\prime$ so that $\phi\alpha=\alpha^\prime \phi$.
  \par\noindent (ii) $[F\oplus F^\prime, \alpha \oplus \alpha^\prime]=
                [F,\alpha]+[F^\prime, \alpha^\prime]$
                \par\noindent (iii) $[F,\alpha\beta]=[F,\alpha]+[F,\beta].$
\end{dfn}
\begin{dfn} We define $K_1(\textbf{Z}\mathcal{P}Hol(X,\mathcal{F})^{bdd})$
  to be the abelian group generated by $[F,\alpha]$ where $F$ is a controlled
  free module and $\alpha$ is an automorphism of $F$ so that
  \par\noindent (i) $[F,\alpha]=[F^\prime, \alpha^\prime]$ if there is an
  isomorphism $\phi:F \to F^\prime$ so that $\phi\alpha=\alpha^\prime \phi$.
  \par\noindent (ii) $[F\oplus F^\prime, \alpha \oplus \alpha^\prime]=
                [F,\alpha]+[F^\prime, \alpha^\prime]$
                \par\noindent (iii) $[F,\alpha\beta]=[F,\alpha]+[F,\beta].$
\end{dfn}
\begin{dfn} $Wh_{\mathcal{F},bg}(\mathcal{P}Hol(X))$ is defined to be the
  quotient of $K_1(\textbf{Z}\mathcal{P}Hol(X)^{bg})$ by the subgroup of
  elements of the form
  $$[F(\sigma), u_{F(\sigma)}]\mbox{ and }[F(\sigma), F(\alpha, \nu)]$$
  where $(S,\alpha)$ is any $bg$ basis over $\mathcal{P}Hol(X,\mathcal{F})$,
  $u_{F(\sigma)}$ is multiplication by a unit and $F(\alpha, \nu)$ is an
  automorphism of bases.
\end{dfn}
\begin{dfn} $Wh_{\mathcal{F},bdd}(\mathcal{P}Hol(X))$ is defined to be the
  quotient of $K_1(\textbf{Z}\mathcal{P}Hol(X)^{bdd})$ by the subgroup of
  elements of the form
  $$[F(\sigma), u_{F(\sigma)}]\mbox{ and }[F(\sigma), F(\alpha, \nu)]$$
  where $(S,\alpha)$ and any $bdd$ basis over $\mathcal{P}Hol(X,\mathcal{F})$,
  $u_{F(\sigma)}$ is multiplication by a unit and $F(\alpha, \nu)$ is an
  automorphism of bases.
  \end{dfn}
\begin{dfn} Let $X, \mathcal{F}$ be a foliated space. Let
  $DR^{bg}_\mathcal{F}(X)$ be the collection of all pairs $(Y,X)$, where
  $Y,\mathcal{G}$ is a foliated space, for which
  $X$ is a foliated $bg$ strong deformation retract of $Y$.
  \end{dfn}

\begin{thm} Every element of $Wh_{\mathcal{F}}(X)$ has a representative of
  the form
  $$(Y,q)=(X,p)\cup_{\phi_r}(S_r \times I^r,q_r)\cup_{\phi_{r+1}}(S_{r+1}
    \times I^{r+1}, q_{r+1}),$$
    for controlled $r$-cell $(S_r \times I^r, q_r)$ and controlled $r+1$-cell
    $(S_{r+1}\times I^{r+1},q_{r+1})$ and
    attaching maps $\phi_r$ and $\phi_{r+1}$.
\end{thm}
\textit{Proof:} The proof is the same as in \cite{Att} Theorem 3.11.
    \begin{thm} The group $Wh^{bg}_{\mathcal{F}}(X)$ is isomorphic to the group
      $Wh_{\mathcal{F}}(\mathcal{P}Hol(X)^{bg})$. The group
      $Wh^{bdd}_{\mathcal{F}}(X)$ is isomorphic to the group
      $Wh_{\mathcal{F}}(\mathcal{P}Hol(X))$.
  \end{thm}
    \textit{Proof:} Observe that $H_i^c$ of the universal cover of a leaf $L$ of
    $\mathcal{F}$ in $DR^{bg}(L)$ is a $\textbf{Z}\mathcal{P}G_1(L)$ module.
    We let the torsion of a chain complex of such modules which is acyclic be
    the class of the boundary map in $Wh_{\mathcal{F}}(\mathcal{P}Hol(X)^{bg})$.
    More specifically, define for $bg$ bases $b_i=((S_i,\sigma_i))$ the
    element $[b_2/b_1]=[F(\sigma_2),F(\alpha,\nu)\phi_1^{-1}\phi_2]$, where
    $(\alpha,\nu)$ is an isomorphism. Now define the torsion of a $bg$
    foliated chain complex: by noting that $B_i=Im \partial_i$ is stably
    free, so using the basis $b_i$,
    $$\tau(C_*)=\sum_i (-1)^i[b_i h_i b_{i-1}/c_i].$$
    This is clearly invariant under leafwise uniform subdivision of $C_*$. Now
    define the isomorphism $\tau:Wh_{\mathcal{F}}(M) \to Wh_{\mathcal{F}}(\mathcal{P}Hol(M)^{bg})$
    by $\tau([X,Y,q])=\tau(i^!C_*^{c,cell}(Y,X))$. This is an isomorphism.

    The same argument works in the bounded case.
    \section{Algebraic L-theory}
    \begin{dfn} Let $(M,\mathcal{F})$ be a manifold $M$ with foliation
      $\mathcal{F}$. Let
      $\mathcal{A}$ be an additive category, $\pi$ a group. Then the
      category $C_M^{\mathcal{F}}(\mathcal{A}[\pi])$ is defined to the
      one whose objects are formal direct sums
      $$M=\sum_{x \in BG}\sum_{y\in L}M(y)$$
      of objects $M(y)$ in $\mathcal{A}[\pi]$, where $L$ is the leaf
      corresponding to the point $x$ in the classifying space of the
      holonomy group $BG$ which fall into a fixed finite
      number of types inside of each ball of fixed radius in $L$ and
      transverse to $L$. Here
      $\mathcal{A}[\pi]$ is the category with the one object $M[\pi]$ for
      each object $M$ in $\mathcal{A}$, and with morphisms linear combinations
      of morphisms $f_g:M\to N$ in $\mathcal{A}$ finite. 
       \end{dfn}
    \begin{dfn} The algebraic surgery group $L^{\mathcal{F},bg}_*(M)$ is defined as
      $$L^{\mathcal{F},bg}_*(M)=L_*(\textbf{Z}\mathcal{P}Hol(M)^{bg})$$
    \end{dfn}
    \begin{dfn} The algebraic surgery group $L^{\mathcal{F},bdd}_*(M)$ is defined as
      $$L^{\mathcal{F},bdd}_*(M)=L_*(\textbf{Z}\mathcal{P}Hol(M)^{bdd})$$
      \end{dfn}
    \begin{dfn} \cite{Ran} A chain duality $(T,e)$ on an additive category $\mathcal{A}$
      is a contravariant functor $T:\mathcal{A}\to \mathcal{B}(\mathcal{A})$
      together with a natural transformation
      $$e:T^2 \to 1:\mathcal{A}\to \mathcal{B}(\mathcal{A})$$
      such that for each object $A$ in $\mathcal{A}$
      \begin{enumerate}
      \item  $e(T(A))\cdot T(e(A))=1:T(A)\to T^3(A)\to T(A).$
      \item  $e(A):T^2(A) \to A$ is a chain equivalence.
      \end{enumerate}
    \end{dfn}
    \begin{dfn} \cite{Ran} Use the standard free $\textbf{Z}[\textbf{Z}_2]$-modules
      resolution of $\textbf{Z}$ to define for any finite chain complex $C$
      in $\mathcal{A}$ the $\textbf{Z}$-module chain complex
      $$W_{\%}C=W\otimes_{\textbf{Z}[\textbf{Z}_2]}(C\otimes_{\mathcal{A}}C)$$
    \end{dfn}
    \begin{dfn} \cite{Ran} An $n$-dimensional quadratic Poincar\'e complex in $\mathcal{A}$
      $(C,\psi)$ is a finite chain complex $C$ in $\mathcal{A}$ together with
      an $n$-cycle $\psi \in (W_\%)_n$ such that the chain map
      $$(1+T)\psi_0:C^{n-*}\to C$$
      is a chain equivalence in $\mathcal{A}$.
    \end{dfn}
    \begin{dfn} \cite{Ran} Let $f:C \to D$ be a chain map. The algebraic mapping cone of
      $f$, $C(f)$ is defined by
      $$d_{C(f)}=\begin{pmatrix} d_D & (-1)^{r-1} f\\ 0 & d_C \end{pmatrix}$$
      $$C(f)_r=D_r\oplus C_{r-1}\to C(f)_{r-1}$$
    \end{dfn}
    \begin{dfn} \cite{Ran} A chain map $f:C\to D$ of finite chain complexes in
      $\mathcal{A}$ induces $\textbf{Z}[\textbf{Z}_2]$-module chain map
      $$f \otimes f:C\otimes_{\mathcal{A}}C\to D \otimes_{\mathcal{A}} D$$
      and hence a $\textbf{Z}$-module chain map
      $$f_{\%}:W_\%C \to W_\%D$$
      \end{dfn}
      \begin{dfn} \cite{Ran} Given an additive category $\mathcal{A}$ and a simplicial complex
      $K$ are combined to define an additive category $\mathcal{A}^*(K)$ of
      $K$-based objects in $\mathcal{A}$ which depends contravariantly on $K$.
      We define $\mathcal{A}_*(K)$ of $K$-based objects in $\mathcal{A}$ which
      depends covariantly on $K$. An object $M$ in an additive category
      $\mathcal{A}$ is $K$-based if it is expressed as a direct sum
      $$M=\sum_{\sigma \in K} M(\sigma)$$
      of objects $M(\sigma)$ in $\mathcal{A}$, such that $\{\sigma \in K \mid
      M(\sigma)\ne 0\}$ is finite. A morphism $f:M\to N$ of $K$-based objects
      is a collection of morphisms in $\mathcal{A}$
      $$f=\{f(\tau,\sigma):M(\sigma)\to N(\tau)\mid \sigma,\tau \in K\}$$
      Let $\mathcal{A}^*(K)$ be the additively category of $K$-based objects
      $M$ in $\mathcal{A}$, with morphisms $f:M\to N$ such that
      $f(\tau, \sigma):M(\sigma)\to N(\tau)$ is 0 unless $\tau \le \sigma$, so
      that
      $$f(M(\sigma))\subseteq \sum_{\tau \le \sigma}N(\tau)$$
      Let $\mathcal{A}_*(K)$ be the additive category of $K$-based objects $M$
      with morphisms $f:M \to N$ so that $f(\tau,\sigma):M(\sigma) \to N(\tau)$ is 0
      unless $\tau \ge \sigma$ so that
      $$f(M(\sigma)\subseteq \sum_{\tau \ge \sigma}N(\tau)$$
      Forgetting the K-based structure defines the covariant assembly functor
      $$\mathcal{A}^*(K)\to\mathcal{A}; M\to M^*(K)=\sum_{\sigma \in K}M(\sigma)$$
      $$\mathcal{A}_*(K)\to \mathcal{A}; M\to M_*(K)=\sum_{\sigma \in K}M(\sigma)$$
      \end{dfn}
    \begin{dfn} \cite{Ran} Let $\mathcal{A}^*[K]$ be the additive category with objects the
      covariant functors
      $$M:K \to \mathcal{A}; \sigma \to M[\sigma]$$
      such that $\{\sigma \in K \mid M[\sigma]\ne 0\}$ is finite. The
      morphisms are the natural transformations of such functors.
      Let $\mathcal{A}_*[K]$ be the additive category with objects the
      contravariant functors as above.
      \end{dfn}
    \begin{dfn}\cite{Ran}
      Let $\mathcal{A}(R)$ be the category of finitely generated free
      $R$-modules, $R$ a ring. The additive category $\mathcal{A}(R)_*[K]$
      will be denoted by $\mathcal{A}[R,K]$ and the additive category
      $\mathcal{A}(R)_*(K)$ will be denoted $\mathcal{A}(R,K)$.
      Their objects will be called $[R,K]$-modules and $(R,K)$-modules.
      \end{dfn}
    \begin{dfn} \cite{Ran} For any additive category with chain duality $\mathcal{A}$ there is
      an algebraic bordism category
      $$\Lambda(\mathcal{A})=(\mathcal{A},\mathcal{B}(\mathcal{A}),\mathcal{C}(\mathcal{A}))$$
      with $\mathcal{B}(\mathcal{A})$ the category of finite chain complexes
      in $\mathcal{A}$, and $\mathcal{C}(\mathcal{A})\subseteq \mathcal{B}(\mathcal{A})$
      the subcategory of contractible complexes.
    \end{dfn}
    \begin{dfn} \cite{Ran} A subcategory $\mathcal{C}\subseteq\mathcal{B}(\mathcal{A})$
        is closed if it is a full additive subcategory which is invariant
        under $T$, such that the algebraic mapping cone $C(f)$ of any chain
        map $f:C\to D$ in $\mathcal{C}$ is an object in $\mathcal{C}$.
\end{dfn}
        \begin{dfn}\cite{Ran}
        A chain complex $C$ in $\mathcal{A}$ is $\mathcal{C}$-contractible if
        it belongs to $\mathcal{C}$. A chain map $f:C \to D$ is
        a $\mathcal{C}$-equivalence
        if the algebraic mapping cone $C(f)$ is $\mathcal{C}$-contractible.
        \end{dfn}
\begin{dfn} \cite{Ran}
        An $n$-dimensional quadratic complex $(C,\psi)$ in $\mathcal{A}$ is
        $\mathcal{C}$-Poincar\'e if the chain complex
        $$\partial C=S^{-1}C((1+T)\psi_0:C^{n-*}\to C)$$
        is $\mathcal{C}$-contractible.
        \end{dfn}
        \begin{dfn} \cite{Ran} Let $\Lambda=(\mathcal{A}, \mathcal{B},\mathcal{C})$ be an
      algebraic bordism category. An $n$-dimensional quadratic complex
      $(C,\psi)$ in $\Lambda$ is an $n$-dimensional quadratic complex in
      $\mathcal{A}$ which is $\mathcal{B}$-contractible and $\mathcal{C}$-Poincar\'e.
      The quadratic $L$-group $L_n(\Lambda)$ is the cobordism group of
      $n$-dimensional quadratic complexes in $\Lambda$.
    \end{dfn}
    \begin{dfn} Let $\mathcal{A}^{\mathcal{F}}_*(K)$ be the additive
      category of $K$-based objects in $\mathcal{A}^{\mathcal{F}}$ with
      morphisms $f:M\to N$ such that $f(\tau,\sigma)=0:M(\sigma)\to N(\tau)$
      unless $\tau \ge \sigma$ so that $f(M(\sigma))\subseteq \sum_{\tau \ge \sigma}N(\tau)$.
    \end{dfn}
    \begin{dfn} The quadratic foliated structure groups of $(R,K)$ are
      the cobordism groups
      $$\mathcal{S}_n^{\mathcal{F}}(R,K)=L_{n-1}(\mathcal{A}^{\mathcal{F}}(R,K),\mathcal{C}^{\mathcal{F}}(R,K),\mathcal{C}^{\mathcal{F}}(R)_*(K))$$
      \end{dfn}
      \begin{dfn} Define the local, uniformly finite, finitely generated free
        $(R,K)$-modules
        $$\Lambda(R)^{\mathcal{F}}_*(K)=(\mathcal{A}^{\mathcal{F}}(R,K),
        \mathcal{B}^{\mathcal{F}}(R,K),\mathcal{C}^{\mathcal{F}}(R)_*(K))$$
        where $\mathcal{B}^{\mathcal{F}}(R,K)$ is the category of finite
        chain complexes of f.g. free foliated $(R,K)$-modules.
        An object in $C^{\mathcal{F}}(R)_*(K)$ is finite f.g. free foliated
        $(R,K)$-module chain complex $C$ such that each $[C][\sigma] (\sigma \in K)$
        is a contractible f.g. free $R$-module chain complex.
      \end{dfn}
      \begin{dfn} Let $M$ be a compact manifold with foliation $\mathcal{F}$.
        Define the foliated normal invariants for $(M,\mathcal{F})$ to be:
        $$H_n(BG;\mathcal{L})=L_n(\Lambda(\textbf{Z})_*^{\mathcal{F}}(M))$$
        where $BG$ is the classifying space of the holonomy groupoid of $M$,
        and $\mathcal{L}$ is the cosheaf assigning to each point $x$ of $BG$
        the $L$-homology of the leaf $S_x$ through $x$.
        \end{dfn}
      \begin{thm}[Surgery Exact Sequence] Let $(M,\mathcal{F})$ be a manifold
        $M$ with a foliation $\mathcal{F}$. Then we have an algebraic surgery
        exact sequence
        $$...\to H_n(BG;\mathcal{L})\to L_n^{\mathcal{F},bg}(M)\to S_n^{\mathcal{F},bg}(R,M)\to
        H_{n-1}(BG;\mathcal{L})\to ...$$
        where $BG$ is the classifying space of the holonomy groupoid and
        $\mathcal{L}$ is the cosheaf of the uniformly finite homology
        with coefficients in the the $L$-spectrum of the leaf $S_x$ through $x\in BG$.
      \end{thm}
      \begin{thm}[Haefliger Cor.3.2.4] \cite{Haefliger} Let $\mathcal{F}$ be a foliation
        on a manifold $X$ so that the holonomy coverings of the leaves are
        $(k-1)$-connected. Then the holonomy groupoid $G$ of $\mathcal{F}$
        considered as a $G$-principal bundle with base $X$ by the end
        projection is $k$-universal. Thus the space $X$
        itself is $k$-classifying and the map $i$ of $X$ to $BG$ is
        $k$-connected.
      \end{thm}
      \begin{thm}[Haefliger Cor.3.1.5] \cite{Haefliger} Suppose the target projection of the
        holonomy groupoid $G$ of the foliation $\mathcal{F}$ on a manifold
        $X$ is a locally trivial fibration, whose fiber is the common
        holonomy covering $L$ of all the leaves. Then the map
        $i:X \to BG$ is homotopy equivalent to a locally trivial fibration
        with base $BG$ and fiber $L$.
      \end{thm}
      The hypothesis of this theorem are satisfied in the following cases:
      \par\noindent i. $X$ is compact and $\mathcal{F}$ possesses a
      transverse Riemannian metric.
      \par\noindent ii. The leaves of $\mathcal{F}$ are transverse to the
      fibers of a compact fibration and the foliation is analytic.
      \par\noindent iii. The leaves of $\mathcal{F}$ are the trajectories
      of a flow without closed orbits or the holonomy group of each closed
      orbit is infinite.
      
      Because of this theorem of Haefliger, foliated surgery has as a special case blocked surgery \cite{Q}, and in analogy with the blocked surgery
      diagram at the end of \cite{Q} we have the following exact sequence:
      \begin{thm}[Blocked surgery exact sequence]
      $$...\to H_n(BG;\mathcal{L})\to H_n(BG;\mathcal{L}^{bg}(S_x)) \to H_n(BG;\mathcal{S}^{bg})\to ...$$
      where $\mathcal{L}$ is cosheaf of the uniformly finite homology with
      coefficients in the $L$ spectrum of the leaf $S_x$ through $x$ in
      $BG$, $\mathcal{L}^{bg}$ is the cosheaf of the $bg$ $L$-theory of the leaf
      $S_x$ through $x$ in $BG$. $\mathcal{S}^{bg}$ is the $bg$ structure
      set of $S_x$ the leaf through $x$ in $BG$.
      \end{thm}
      In addition there are exact sequences
\begin{thm}[Leafwise Assembly]
      $$...\to Fiber^{NI}_n \to H_n(BG;\mathcal{L}) \overset{A^{NI}}{\to} H_n(M;\textbf{L})\to Fiber^{NI}_{n-1}\to ...$$
      $$...\to Fiber^{L}_n \to H_n(BG;\mathcal{L}^{bg}) \overset{A^{L}}{\to} L_n(\pi_1(M)) \to Fiber^{L}_{n-1} \to ...$$
      $$... \to Fiber^{S}_n \to H_n(BG;\mathcal{S}^{bg}(S_x))\overset{A^{S}}{\to} \mathcal{S}_n(M) \to Fiber^{S}_{n-1} \to ...$$
      where $Fiber^{NI}_n, Fiber^{L}_n, Fiber^{S}_n$ are the fibers of
      the leafwise assembly maps, and classify the non-leaves,
      and the maps $A^{NI}:H_n(BG;\mathcal{L}) \to H_n(M;\textbf{L})$,
      $A^{L}:H_n(BG;\mathcal{L}^{bg}) \to L_n(\pi_1(M))$
      and $A^S:H_n(BG;\mathcal{S}^{bg}(S_x)) \to \mathcal{S}_n(M)$ are the leafwise
      assembly maps.
\end{thm}
Note that the assembly map $A^{NI}:H_n(BG;\mathcal{L})\to H_n(M;\textbf{L})$
was used in \cite{AttCap1} to prove Bott Integrability in the simply-connected
case. It was generalized in \cite{AttCap1} and \cite{AttCap2} to prove more
general versions of Bott Integrability.
\begin{thm}[Embeddability of Leaves]
      The image of $Fiber_n^{S} \to H_n(BG;\mathcal{S}^{bg}(S_x))$
      is the kernel of $A^{S}:H_n(BG;\mathcal{S}^{bg}(S_x))\to \mathcal{S}_n(M)$
      and therefore we take \break $H_n(BG;\mathcal{S}^{bg})/Ker A^{S}$ to be the
      foliations on $M$ with the leaves of a given $bg$ homotopy type. This
      answers to a certain extent the question of when an open manifold embeds
      as a leaf of a foliation in a given closed manifold $M$.
\end{thm}
      For example,
      take $N \times \mathcal{F}$, where $\mathcal{F}$ is the Kronecker
      foliation on $T^2$, a foliation on $N \times T^2$. $BG=T^2$, the leaves
      are $bg$ homotopy equivalent to $N \times \textbf{R}$,
      $$\mathcal{S}^{bg}(N \times \textbf{R})=H_*^{uf}(\textbf{R};\mathcal{S}_*(N))$$
      and $Ker A^{S}$ is the set of non-end-periodic leaves \cite{AH}. This gives us
      the Pimsner-Voiculescu theorem from Section 2, as the set of end-periodic
      leaves corresponds to $\textbf{Z} \subset H_0^{uf}(\textbf{R})$ as in
      \cite{AH} and therefore the structure set is $H_*(BG;\mathcal{S}(N))=H_*(T^2;\mathcal{S}(N))$.
      \begin{thm}[Bounded Surgery Exact Sequence]
        Let $(M,\mathcal{F})$ be a manifold $M$ with a foliation $\mathcal{F}$.
        Then we have an algebraic surgery exact sequence
        $$...\to H_n(BG;\mathcal{L})\to L_n^{\mathcal{F},bdd}(M)\to S_n^{\mathcal{F},bdd}(R,M)\to
          H_{n-1}(BG;\mathcal{L})\to ...$$
          where $BG$ is the classifying space of the holonomy groupoid and
          $\mathcal{L}$ is the cosheaf of the locally finite homology
          of the $L$-spectrum of the leaf $S_x$ through $x \in BG$.
          \end{thm}
      \newtheorem{conj}{Conjecture}[section]
      We have two Novikov and Borel Conjectures for foliations. The first,
      ones for bounded surgery, correspond to the Baum-Connes Conjecture
      for foliations and the Novikov Conjecture for foliations due to
      Baum and Connes \cite{BC}.
      \begin{conj}[Bounded Novikov Conjecture for Foliations]
        The map
        $$H_*(BG;\mathcal{L}^{lf})\to L_*^{\mathcal{F},bdd}(M)$$
        where $\mathcal{L}^{lf}$ is the cosheaf assigning the locally finite $L$-homology
        of the leaf through $x \in BG$ to the point $x \in BG$,
        is injective, provided the leaves of $(M,\mathcal{F})$ are
        uniformly contractible.
      \end{conj}
      \begin{conj}[Bounded Borel Conjecture for Foliations]
        If the leaves of $(M,\mathcal{F})$ are uniformly contractible, then
        the map
        $$H_*(BG;\mathcal{L}^{lf})\to L_*^{\mathcal{F},bdd}(M)$$
        where $\mathcal{L}^{lf}$ is the cosheaf assigning the locally finite $L$-homology
        of the leaf through $x \in BG$ to the point $x \in BG$
        is an isomorphism.
      \end{conj}
      The second two conjectures are for bounded geometry surgery:
      \begin{conj}[Bounded Geometry Novikov Conjecture for Foliations]
        If the leaves of $(M,\mathcal{F})$ are uniformly contractible, then
        $$H_*(BG;\mathcal{L}^{uf}(S_x))\to L_*^{\mathcal{F},bg}(M)$$
        is injective, where $S_x$ is the leaf of $\mathcal{F}$ through $x \in BG$
        and $\mathcal{L}^{uf}$ is the uniformly finite $L$-homology cosheaf
        of the leaves over $BG$.
      \end{conj}
      \begin{conj}[Bounded Geometry Borel Conjecture for Foliations]
        If the leaves of $(M,\mathcal{F})$ are uniformly contractible, then
        $$H_*(BG;\mathcal{L}^{uf}(S_x))\to L_*^{\mathcal{F},bg}(M)$$
        is an isomorphism, where $S_x$ is the leaf of $\mathcal{F}$ through
        $x \in BG$ and $\mathcal{L}^{uf}$ is the uniformly finite $L$-homology
        cosheaf of the leaves over $BG$.
      \end{conj}
        \section{Examples of Foliations Satisfying the Novikov and Borel Conjectures}
      We recall the cases where the Baum-Connes conjecture has been verified.
We will prove the bounded Borel conjecture for these foliations.
      \par\noindent 1. Fibrations $F \to M \to B$. In this case the leaf space
      is identified with the base space of the fibration $B$.
      \par\noindent 2. Foliations induced by free actions of $\textbf{R}^n$,
      and free actions of a solvable simply connected Lie group $\Gamma$.
      \par\noindent 3. The Reeb foliations on $T^2$ and $S^3$.
      \par\noindent 4. Foliations without holonomy. In this case $BG$ is a
      torus $T^n$ provided with a linear foliation, and we apply 2.
      \par\noindent 5. Foliations almost without holonomy. Applying graphs of
      groups, this is reduced to 4.
      \par\noindent 6. Nontrivial examples: the Sacksteder foliation, the
      Hirsch foliation, $\textbf{Z}$-periodic foliations.
\bigskip\par\noindent
      Our Pimsner-Voiculescu theorem is a case of the Borel Conjecture for
      foliations, the case of the Kronecker foliation. We will verify the
      Borel Conjecture in all cases where the Baum-Connes Conjecture is known.
      \begin{thm} The Borel Conjecture is true for fibrations
        $$F\to B\to M$$
        where $F$ is compact.
      \end{thm}
      \textit{Proof:} In this case the leaf space is $B$. We see that
      $$H_n(BG;\underline{H_m^{lf}(F;\textbf{L})})=H_n(B;\underline{H_m^{lf}(F;\textbf{L})})=H_{n+m}(M;\textbf{L})$$
      Where $\underline{H_*^{lf}(F;\textbf{L})}$ is the constant cosheaf.
      The surgery group of the foliation is isomorphic to this.

      \begin{thm} The Borel Conjecture is true for the Reeb foliation on
        $T^2$ and $S^3$
      \end{thm}
      \textit{Proof:} This can be seen easily by codimension 1 splitting
      as for the case of Pimsner-Voiculescu.
      \begin{thm} The Borel Conjecture is true for free actions of $\textbf{R}^n$
        and free actions of a solvmanifold $\Gamma$.
      \end{thm}
      \textit{Proof:} We again use codimension 1 splitting along the leaves
      inductively down to $\textbf{R}^{n-1}$ etc.
      \begin{thm} The Borel Conjecture is true for foliations without holonomy.
      \end{thm}
      \textit{Proof:} Codimension 1 splitting applies where $BG$ is $T^n$ and
      the leaves are linear.
      \begin{thm} The Borel Conjecture is true for foliations almost without
        holonomy.
      \end{thm}
      \textit{Proof:} Apply the previous using graphs of groups.
      \begin{thm} The Borel Conjecture is true for the Hirsch foliation.
\end{thm}
        \textit{Proof:} Apply codimension 1 splitting to the leaves, which
        are in the form of an $n$-partite graph \cite{BisHurderShive}.

        Similarly, the Borel conjecture can be verified for the Sacksteder
        foliation and $\textbf{Z}$-periodic foliations from codimension 1
        splitting.

        \bigskip\par\noindent
The Novikov conjecture for foliations has been stated by P.Baum and A.Connes \cite{BC}. Let $\pi$ be the fundamental groupoid along the leaves. $B\pi$ denotes
the classifying space of the topological groupoid $\pi$. Since $V$ is the
units of $\pi$ there is a canonical map:
$$\lambda: V \to B\pi$$
$\pi$ is itself a principal $\pi$-bundle over $V$ and the map $\lambda$
is the classifying map of this principal $\pi$-bundle.
\begin{conj} Let $(V,F)$ and $(V^\prime,F)$ be orientable $C^\infty$
  foliations with $V,V^\prime$ compact. Let $f:V^\prime \to V$ be a
  leafwise homotopy equivalence. Choose orientations for $V$ and $V^\prime$
  so that $f$ is orientation preserving. Then in $H_*(B\pi;\textbf{Q})$ there
  is the equality:
  $$\lambda_*(\textbf{L}(TV)\cap [V])=\lambda_* f_*(\textbf{L}(TV^\prime)\cap [V^\prime])$$
  \end{conj}
Baum and Connes proved the Foliated Novikov Conjecture in the following special case:
      \par\noindent Let $\mathcal{L}$ be a leaf of the foliation $(V,F)$.
      If for every leaf $\mathcal{L}$,
      $\pi_i(\mathcal{L})=0$ for all $i\ge 2$, then $B\pi=V$, where $B\pi$ is
      the classifying space of the fundamental  groupoid along the leaves.
      If $V$ is compact
      and $(V,F)$ has negatively curved leaves this is the case, so the
      Novikov conjecture holds.

      S.Hurder \cite{Hurder}
      has extended this to the class of ultra-spherical foliations,
      which will be defined below.
      
      \begin{dfn} \cite{Hurder} Let $(M,\mathcal{F})$ be a foliation and
        $\mathcal{G}_{\mathcal{F}}$ be its holonomy groupoid.
        A leafwise path $\gamma$ is a continuous
        map $\gamma:[0,1]\to M$ whose image is contained in a single leaf
        of $\mathcal{F}$. There are natural continuous maps
        $s,r:\mathcal{G}_{\mathcal{F}}\to M$ defined by $s(\gamma)=\gamma(0)$
        and $r(\gamma)=\gamma(1)$. For a point $x\in M$, the pre-image
        $s^{-1}(x)=\tilde{L}_x$ is the holonomy cover of the leaf
        $L_x$ of $\mathcal{F}$ through $x$.
        \end{dfn}
      \begin{dfn} \cite{Hurder} For $x \in M$ and a leafwise path $\gamma:[0,1] \to \tilde{L}_x$,
        define the plaque length function $\mathcal{N}_T(\gamma)$ to be the
        least number of plaques to cover the image of $\gamma$. Define the
        plaque distance function $D_x(\cdot,\cdot)$ on the holonomy cover
        $\tilde{L}_x$ using the plaque length function: for $y,y^\prime \in
        \tilde{L}_x$
        $$D_x(y,y^\prime)=\inf\{\mathcal{N}_T \mid \gamma \mbox{ is a leafwise
          path from y to }y^\prime\}$$
      \end{dfn}
      \begin{dfn} \cite{Hurder} Denote by $C_u(\mathcal{G}_{\mathcal{F}})\subset C(\mathcal{G}_{\mathcal{F}})$
        the subspace of uniformly continuous functions on $\mathcal{G}_{\mathcal{F}}$,
        the holonomy groupoid of $\mathcal{F}$.
        \end{dfn}
      \begin{dfn} \cite{Hurder} For $x\in M$ and $r>0$, define the fiberwise variation function
        $$V_s(x,r):C(\tilde{L}_x) \to [0,\infty]$$
        $$V_s(x,r)(h)(y)=sup\{\mid h(y^\prime)-h(y)\mid\mbox{ such that }
        D_x(y,y^\prime)\le r\}$$
        We say that $f \in C(\mathcal{G}_{\mathcal{F}})$ has uniformly
        vanishing variation at infinity, where $\mathcal{G}_{\mathcal{F}}$
        is the holonomy groupid of $\mathcal{F}$ if there exists a function
        $D:(0,\infty) \to [0,\infty)$ so that if $D_x(y,*x) > D(\epsilon)$
          then $V_s(x,r)(i_x^*f)(y) < \epsilon$. Here $*x$ denotes the
          constant path at $x$. Let $C_h(\mathcal{F})
          \subset C_u(\mathcal{F})$ denote the subspace of uniformly
          continuous functions which have uniformly vansishing variation
          at infinity.
      \end{dfn}
      
      \begin{dfn} \cite{Hurder} Let $\mathcal{F}$ be a topological foliation of a paracompact
        manifold $M$ equipped with a regular foliation atlas. The corona,
        $\partial_h\mathcal{F}$, of $\mathcal{F}$ is the spectrum of the
        quotient $C^*$-algebra $C_h(\mathcal{F})/C_0(\mathcal{F})$.
        \end{dfn}
      \begin{dfn} \cite{Hurder} A foliation $\mathcal{F}$ on a connected manifold $M$ is said
        to be ultra-spherical if there exists a map of coronas
        $\sigma:\partial_h \mathcal{F} \to S\mathcal{F}$ which commutes with
        the projections onto $M$, and so the $\sigma^*\Theta \in H^*(\partial_h
        \mathcal{F})$ is nonzero.
        \end{dfn}
      We next prove the foliated Novikov Conjecture using bounded surgery.
      \begin{thm} Let $\mathcal{F}$ be an oriented ultra-spherical foliation
        with uniformly contractible leaves and Hausdorff holonomy groupoid.
        Then the foliated Novikov Conjecture is true for $\mathcal{F}$.
      \end{thm}
      \textit{Proof:}See \cite{Hurder}.
      Let
      $$\mu:KO(BG)\otimes \textbf{Q} \to L^{\mathcal{F},bdd}_*(M)\otimes \textbf{Q}$$
      be the assembly map, which we must show to be injective. There
      exists a boundary class $u \in KO^*(\partial_h \mathcal{F})$.
      Let $\eta \in KO(S\mathcal{F})$ with KO-theory boundary
      $\beta[T\mathcal{F}] \in KO(T\mathcal{F})$ and set
      $[u]=\sigma^*\eta$. There is a continuous extension of $\sigma$ to
      a map of pairs
      $$\overline{\sigma}:(\overline{\mathcal{G}_{\mathcal{F}}},\partial_h \mathcal{F})
      \to (\overline{T\mathcal{F}}, S\mathcal{F})$$
      which commutes with the projection onto $M$. By naturality of the
      boundary map, $\partial [u]=\overline{\sigma}^*\beta[T\mathcal{F}]$.
      The rest follows as in \cite{Roe}, Theorem 6.9.
      \section{Further Directions}
      The difference between boundedly controlled foliated surgery and bounded
      geometry foliated surgery should give rise to a fiber (to the
      forgetful map) and an exact sequence. The fiber is the set of
      metrics on the leaves given the boundedly controlled topological type.
      We have not investigated the situation where a foliation has a
      symmetry or group action, and the resulting equivariant surgery theory.
      Furthermore, the symmetric L-groups due to Mishchenko and Ranicki
      \cite{Ran}, and visible L-groups due to Weiss \cite{Ran} also merit
      investigation.

\end{document}